
\input amstex
\documentstyle{amsppt}
\magnification1200 
\pagewidth{6.5 true in} 
\pageheight{9.25 true in}
\NoBlackBoxes
\def\Lam{\Lambda}

\topmatter
\title Pretentious multiplicative functions and an inequality 
for the zeta-function \endtitle
\author Andrew Granville and K. Soundararajan
\endauthor
\rightheadtext {Pretentious multiplicative functions}

\address{D{\'e}partment  de Math{\'e}matiques et Statistique,
Universit{\'e} de Montr{\'e}al, CP 6128 succ Centre-Ville,
Montr{\'e}al, QC  H3C 3J7, Canada}\endaddress
\email{andrew{\@}dms.umontreal.ca}
\endemail
\address{Department of Mathematics, University of Michigan, Ann Arbor,
Michigan 48109, USA} \endaddress \email{ksound{\@}umich.edu}
\endemail
\curraddr {Department of Mathematics, Stanford University, 
450 Serra Mall, Bldg. 380, Stanford, CA 94305-2125, USA} 
\endcurraddr
\email{ksound{\@}math.stanford.edu}; ksound{\@}umich.edu\endemail
\thanks{Le premier auteur est partiellement soutenu par une bourse
de la Conseil  de recherches en sciences naturelles et en g\' enie
du Canada. The second  author is partially supported by the American Institute 
of Mathematics and the National Science Foundation.}
\endthanks
\abstract{We note how several central results in multiplicative number theory
may be rephrased naturally in terms of multiplicative functions $f$ that 
pretend to be another multiplicative function $g$.  We formalize a `distance' 
which gives a measure of such {\sl pretentiousness}, and as one consequence obtain 
a curious inequality for the zeta-function.}
\endabstract
\endtopmatter

\def\L{\fracwithdelims()}
\document

\noindent A common theme in several problems in multiplicative number theory 
involves identifying multiplicative functions $f$ that pretend to be another multiplicative 
function $g$.   Indeed, this theme may be found as early as in the proof of the 
prime number theorem; in particular in showing that $\zeta(1+it) \neq 0$.  
For, if $\zeta(1+it)$ equals zero, then we expect the Euler product 
$\prod_{p \le P} (1-1/p^{1+it})^{-1}$ to be small.  This means that $p^{-it} 
\approx -1$ for many small primes $p$; or equivalently, that the multiplicative 
function $n^{-it}$ pretends to be the multiplicative function $(-1)^{\Omega(n)}$.   The insight 
of Hadamard and de la Vallee Poussin is that in such a case $n^{-2it}$ would 
pretend to be the multiplicative function that is identically $1$, and 
this possibility can be eliminated by noting that $\zeta(1+2it)$ is regular for $t\neq 0$.  
 
Another example is given by Vinogradov's conjecture that the least quadratic 
non-residue $\pmod p$ is $\ll p^{\epsilon}$.  If this were false, then the 
Legendre symbol $\L{n}{p}$ would pretend to be the trivial character for a long range  
of $n$.   Even more extreme is the possibility that a quadratic Dirichlet $L$-function 
has a Landau-Siegel zero (a real zero close to $1$), in which case that 
quadratic character $\chi$ would pretend to be the function $(-1)^{\Omega(n)}$.  In 
both these examples, it is not known how to eliminate the possibility of such pretentious behavior 
by characters. 

A third class of examples is provided by the theory of mean values of 
multiplicative functions.  Let $f(n)$ be a multiplicative function
with $|f(n)|\leq 1$ for all $n$, and consider when the mean value
$$
\frac{1}{x} \sum_{n\le x} f(n)\tag{1}
$$
can be large in absolute value; for example, when is it $\gg 1$?
If we write $f(n)=\sum_{d|n} g(d)$ for a multiplicative function $g$, exchange sums, 
and ignore error terms, then we are led to expect that the mean value in (1) is 
about
$$
\prod_{p\leq x} \Big( 1-\frac 1p\Big)
\Big( 1+\frac {f(p)}p+\frac{f(p^2)}{p^2}+\dots \Big),
$$
which has size about 
$$
\exp\Big( -\sum_{p\leq x}   \frac {1-f(p)}p\Big) . \tag{2}
$$
The quantity in (2) is large if and only if $f(p)$ is roughly equal to 1, for ``most'' primes $p\leq x$.  Therefore we may guess that (1) is large only if $f$ pretends to be the constant function $1$.  

When $f$ is non-negative (so $0\le f(n)\le 1$), a result of R.R. Hall [8] gives that (1) is 
$\ll$ (2), confirming our guess.   If we restrict ourselves to real valued $f$ (so $-1\le f(n) \le 1$) 
then another result of Hall [7] gives that 
$$
\frac{1}{x} \sum_{n\le x} f(n) \ll \exp\Big( -\kappa\sum_{p\leq x}   \frac {1-f(p)}p\Big).  
$$
Here $\kappa=0.3286 \ldots$ is an explicitly given constant, 
and the result is false for any larger value of $\kappa$.   
Thus our heuristic that (1) is of size at most (2) does not hold, but nonetheless
our guess that (1) is large only if $f$ pretends to be $1$ is correct.

When $f$ is allowed to be complex valued, another possibilities for (1) being large arises. 
Note that 
$$
\frac{1}{x} \sum_{n\le x} n^{i\alpha}
\sim \frac{x^{i\alpha}}{1+i\alpha},
$$
so that (1) is large in absolute value when $f(n)=n^{i\alpha}$. 
G. Hal{\' a}sz ([5], [6]) made the beautiful realization that this is essentially the 
only way for (1) to be large: that is $f$ must pretend to be the function $n^{i\alpha}$ 
for some real number $\alpha$.    After incorporating significant refinements by Montgomery and 
Tenenbaum, a version of Hal{\' a}sz'a result (see [9]) is that if
$$
M(x,T) := \min_{|t| \le 2T} \sum_{p\le x}
\frac{1-\text{Re} (f(p) p^{-it})}{p} 
$$
then
$$
\frac{1}{x}\Big|\sum_{n\le x} f(n) \Big| \ll (1+M(x,T)) e^{-M(x,T)}
+ \frac{1}{\sqrt{T}}.
$$
For an explicit version of this see [4].

Recently, in [1] A. Balog and the authors considered the mean value of 
multiplicative functions along arithmetic progressions:  that is, for $q < x$ and 
$(a,q)=1$,
$$
\frac{q}{x} \sum\Sb n\leq x \\ n\equiv a \pmod q\endSb f(n). \tag{3}
$$
If $f$ is a character $\chi\pmod q$ then the above is essentially $f(n)=\chi(a)$ for every term in the sum in (3), 
and so the mean value is large.   If we take $f(n)=\chi(n) n^{i\alpha}$ for 
a fixed real number $\alpha$, then also we would get a large mean value.  
In [1] we show, generalizing Halasz's results, that if $q\le x^{\epsilon}$ then 
these are the only ways of getting a 
large mean value in (3).  

These examples suggest that one should   define a distance between multiplicative 
functions, which would quantify how well $f$ pretends to be another function $g$.   We 
formulated such a notion in our recent work on the P{\' o}lya-Vinogradov inequality [3].  
This states (see [2] for example) that for a primitive character $\chi \pmod q$
$$
\max_{x} \Big|\sum_{n\le x} \chi(n)\Big|  \ll \sqrt{q} \log q,  \tag{4}
$$
and in [3] we showed that (4) can be substantially improved unless $\chi$ pretends 
to be a character of much smaller conductor.   The precise characterization 
in fact enabled us to improve (4) in many circumstances, for instance for cubic characters $\chi$.  
In this article we draw attention to this notion of distance, and record some amusing 
inequalities that it leads to.



Consider the space ${\Bbb U}^{{\Bbb N}}$ of vectors 
${\bold z} =(z_1,z_2,\ldots)$ where each $z_i$ lies on the unit disc ${\Bbb U} 
= \{|z|\le 1\}$.  The space is equipped with a product obtained by multiplying componentwise: that is, 
${\bold z}\times {\bold w}=(z_1w_1,z_2w_2,\ldots)$.  Suppose we have a sequence of functions $\eta_j: {\Bbb U} \to {\Bbb R}_{\ge 0}$ for which $\eta_j(zw) \le \eta_j(z) +\eta_j(w)$ for any $z$, $w \in {\Bbb U}$.  
Then we may define a `norm' on ${\Bbb U}^{\Bbb N}$ by setting 
$$
\| {\bold z}\| = \Big(\sum_{j=1}^{\infty} \eta_j(z_j)^2 \Big)^{\frac 12},
$$
assuming that the sum converges.   The key point is that such 
a norm satisfies the triangle inequality 
$$
\| {\bold z} \times {\bold w} \| \le \|{\bold z}\| +\|{\bold w}\|. \tag{5}
$$
Indeed we have 
$$
\align
\| {\bold z}\times {\bold w}\|^2
&= \sum_{j=1}^{\infty}  \eta_j(z_j w_j)^2 
\le \sum_{j=1}^{\infty} (\eta_j(z_j)^2 +\eta_j(w_j)^2 +2 \eta_j(z_j)\eta_j(w_j))
\\
&\le \| {\bold z}\|^2 + \|{\bold w}\|^2 + 2 \Big(\sum_{j=1}^{\infty} \eta_j(z_j)^2 \Big)^{\frac 12} 
\Big( \sum_{j=1}^{\infty} \eta_j(w_j)^2\Big)^{\frac 12} 
= (\|{\bold z}\|+\|{\bold w}\|)^2,\\
\endalign
$$
using the Cauchy-Schwarz inequality, which implies (5).

A nice class of examples is provided by taking $\eta_j(z)^2 = a_j (1-\text{Re }z)$ 
where the $a_j$ are non-negative constants with $\sum_{j=1}^{\infty} a_j <\infty$.  
This last condition ensures the convergence of the sum in the definition of 
the norm.  To verify that $\eta_j(zw) \le \eta_j(z)+\eta_j(w)$, 
note that $1-\text{Re }(e^{2i\pi\theta})=2\sin^2(\pi \theta)$ and   $|\sin(\pi (\theta+\phi))|\leq |\sin(\pi \theta)\cos(\pi \phi)|+|\sin(\pi \phi)\cos(\pi \theta)|\leq |\sin(\pi \theta)|+|\sin(\pi \phi)|$.  
This settles the case where $|z|=|w|=1$, and one can extend this to all pairs 
$z,w\in \Bbb U$.

 Now we show how to use such norms to study multiplicative functions.    Let 
 $f$ be a completely multiplicative function.   Let $q_1< q_2 <\ldots$ denote 
 the sequence of prime powers, and we identify $f$ with the element in 
 ${\Bbb U}^{{\Bbb N}}$ given by $(f(q_1), f(q_2), \ldots)$.  Take 
 $a_j= \Lambda(q_j)/(q_j^{\sigma}\log q_j)$ for $\sigma>1$, and $\eta_j(z)^2=
 a_j(1-\text{Re }z)$.    Then our norm is 
 $$
 \|f\|^2 = \sum_{j=1}^{\infty} \frac{\Lam(q_j)}{q_j^{\sigma}\log q_j} (1-\text{Re }f(q_j)) 
 = \log \frac{\zeta(\sigma)}{|F(\sigma)|},
 $$
 where $F(s)=\sum_{n=1}^{\infty } f(n)n^{-s}$. 
   
\proclaim{Proposition 1}  Let $f$ and $g$ be completely multiplicative functions 
with $|f(n)|\le 1$ and $|g(n)|\le 1$.   Let $s$ be a complex number with Re $s >1$, and 
set $F(s)=\sum_{n=1}^{\infty} f(n)n^{-s}$, $G(s)=
\sum_{n=1}^{\infty} g(n)n^{-s}$, and $F\otimes G(s)= \sum_{n=1}^{\infty} f(n)g(n)n^{-s}$.  
Then, for $\sigma>1$,
$$
\sqrt{\log \frac{\zeta(\sigma)}{|F(\sigma)|}}+\sqrt{\log \frac{\zeta(\sigma)}{|G(\sigma)|}}
\ge \sqrt{\log \frac{\zeta(\sigma)}{|F\otimes G(\sigma)|}},
$$
and 
$$
\sqrt{\log |\zeta(\sigma)F(\sigma)|} +\sqrt{\log |\zeta(\sigma)G(\sigma)|} 
\ge \sqrt{\log \frac{\zeta(\sigma)}{|F\otimes G(\sigma)|}}.
$$
\endproclaim 
\demo{Proof}  The first inequality follows at once from the triangle inequality.  
The second inequality follows upon taking $(-1)^{\Omega(n)}f(n)$ and $(-1)^{\Omega(n)}g(n)$ 
in place of $f$ and $g$, and using the first inequality.  
\enddemo

If we take $f(n)=n^{-it_1}$ and $g(n)=n^{-it_2}$ then we are 
led to the following curious inequalities for the zeta-function which 
we have not seen before. 

\proclaim{Corollary 2} We have
$$
\sqrt{\log \frac{\zeta(\sigma)}{|\zeta(\sigma+it_1)|}}+\sqrt{\log \frac{\zeta(\sigma)}{|\zeta(\sigma+it_2)|}}
\ge \sqrt{\log \frac{\zeta(\sigma)}{|\zeta(\sigma+it_1+it_2)|}},
$$
and 
$$
\sqrt{\log |\zeta(\sigma)\zeta(\sigma+it_1)|} +\sqrt{\log |\zeta(\sigma)\zeta(\sigma+it_2)|} 
\ge \sqrt{\log \frac{\zeta(\sigma)}{|\zeta(\sigma+it_1+it_2)|}}.
$$
\endproclaim

If we take $t_1=t_2$ in the second inequality of Corollary 2, square out and simplify, 
we obtain the classical inequality $\zeta(\sigma)^3 |\zeta(\sigma+it)|^4 |\zeta(\sigma+2it)| 
\ge 1$.  It is conceivable that the more flexible inequalities in Corollary 2 could 
lead to numerically better zero-free regions for $\zeta(s)$, but our initial approaches 
in this direction were unsuccessful.  

Taking $f(n)=\chi(n)n^{-it_1}$ and $g(n)=\psi(n)n^{-it_2}$ in Proposition 1 leads 
to similar inequalities for Dirichlet $L$-functions: for example,
$$
\sqrt{ \log  \frac{ \zeta(\sigma)}{|L(\sigma+it_1+it_2,\chi \psi)|}  }\leq 
\sqrt{ \log \frac{  \zeta(\sigma)}{|L(\sigma+it_1,\chi)|}   }+ \sqrt{\log \frac{ \zeta(\sigma)}{|L(\sigma+it_2,\psi)|}}. 
$$
Thus the classical inequalities leading to zero-free regions for Dirichlet $L$-functions can 
be put in this framework of triangle inequalities.  We wonder if similar useful inequalities 
could be found for other $L$-functions.  
 
It is no more difficult to conclude in Proposition 1 that
$$
\sqrt{ \pm \text{Re} \Big( \frac{F'(\sigma)}{F(\sigma)} \Big) -  \frac{\zeta'(\sigma)}{\zeta(\sigma)} } +
\sqrt{   \pm \text{Re} \Big( \frac{G'(\sigma)}{G(\sigma)} \Big) -\frac{\zeta'(\sigma)}{\zeta(\sigma)}} \ge
\sqrt{ \text{Re} \Big( \frac{(F\otimes G)'(\sigma)}{(F\otimes G)(\sigma)} \Big) - \frac{\zeta'(\sigma)}{\zeta(\sigma)}} .
$$
Again taking $F=G$ and squaring we obtain:
$$
3   \frac{\zeta'(\sigma)}{\zeta(\sigma)}   
\pm 4 \text{Re} \left( \frac{F'(\sigma)}{F(\sigma)} \right)   +\text{Re} \left( \frac{(F\otimes F)'(\sigma)}{(F\otimes F)(\sigma)} \right) \le 0 .
$$

Above we saw one way of defining a norm on multiplicative functions.  Another way 
is to define the distance (up to $x$) between the multiplicative functions $f$ and $g$ by
$$
{\Bbb D}(f,g;x)^2 = \sum_{p\le x} \frac{1-\text{Re }f(p)\overline{g(p)}}{p}.
$$
This arises by taking $a_j=1/q_j$ if $q_j$ is a prime $\le x$, and $a_j=0$ otherwise.  
Thus we have the triangle inequality 
$$
{\Bbb D}(1,f;x) +{\Bbb D}(1,g;x) \ge {\Bbb D}(1,fg;x),
$$
where $1$ denotes the multiplicative function that is $1$ on all natural numbers.  
Notice that this distance came up naturally in our discussion of the results of Hall and
 Hal{\' a}sz on mean values of multiplicative functions.  This distance also 
 provided a convenient framework for our work in [3], where we established 
 the following lower bounds for the distance between characters.  

 \proclaim{Lemma 3} Let $\chi \pmod q$ be a primitive character 
of odd order $g$.  Suppose $\xi \pmod m$ is a primitive character 
such that $\chi(-1)\xi(-1)=-1$.  If $m\le (\log y)^A$ then 
$$ 
{\Bbb D}(\chi,\xi;y)^2 \ge  \left( 1- \frac{g}{\pi} \sin \frac{\pi}{g} +o(1)\right) \log \log y.
$$
\endproclaim

\demo{Proof}  See Lemma 3.2 of [3].
\enddemo

\proclaim{Lemma 4} Let $g\ge 2$ be fixed.  
Suppose that for $1\le j\le g$,  
$\chi_j \pmod{q_j}$ is a primitive character.    
Let $y$ be large, and suppose $\xi_j \pmod {m_j}$ are 
primitive characters with conductors $m_j \le \log y$. 
Suppose that $\chi_1\cdots \chi_g$ is the trivial character, 
but $\xi_1\cdots \xi_g$ is not trivial.  Then 
$$
\sum_{j=1}^{g} {\Bbb D}(\chi_j,\xi_j;y)^2 
\ge \Big(\frac{1}{g} +o(1)\Big) \log \log y.
$$
\endproclaim
 
\demo{Proof}  See Lemma 3.3 of [3].
\enddemo

\proclaim{Lemma 5} Let $\chi \pmod q$ be a primitive 
character.  Of all primitive characters with conductor below 
$\log y$, suppose that $\psi_j \pmod {m_j}$ ($1\le j\le A$) 
give the smallest distances ${\Bbb D}(\chi,\psi_j;y)$ 
arranged in ascending order.  Then for each $1\le j \le A$ 
we have that 
$$
{\Bbb D}(\chi,\psi_j;y)^2 \ge \Big(1-\frac{1}{\sqrt{j}} +o(1)
\Big) \log \log y.
$$
\endproclaim

\demo{Proof}  See Lemma 3.4 of [3].
\enddemo

We conclude this  article by showing, in a suitable sense, that a multiplicative function $f$ cannot 
pretend to be two different characters.  This is in some ways a generalization of 
the fact that there is ``at most one Landau-Siegel zero," which may be viewed as saying that 
$\mu(n)$ cannot pretend to be two different characters
with commensurate conductors.

 \proclaim{Proposition 6} Let $\chi \pmod q$ be a primitive character.  There is an absolute 
 constant $c>0$ such that for all 
 $x\ge q$ we have
 $$
 {\Bbb D}(1,\chi;x)^2 \ge \frac12 \log \Big( \frac{c\log x}{\log q} \Big) .
 $$
Consequently, if $f$ is a multiplicative function, and $\chi$ and $\psi$ are any 
two distinct primitive characters with conductor below $Q$,  then for $x \ge Q$ we 
have
$$
{\Bbb D}(f,\chi;x)^2 + {\Bbb D}(f,\psi;x)^2 \ge \frac{1}{8} \log \Big( \frac{c\log x}{2\log Q} \Big) .
$$
\endproclaim
\demo{Proof} Let $d_{\chi}(n)=\sum_{ab=n} \chi(a)\overline{\chi(b)}$.  Thus 
$d_{\chi}(n)$ is a real valued multiplicative function which satisfies $|d_{\chi}(n)|
\le d(n)$ for all $n$.  We begin by noting that 
$$
\sum_{n\le x} d_{\chi}(n) \ll \sqrt{qx}\log q+ q (\log q)^2. \tag{6}
$$
To prove (6) note that if $n=ab \le x$ then either $a \le \sqrt{x}$ or $b\le\sqrt{x}$ or both.  
Therefore 
$$
\sum_{n\le x} d_{\chi}(n) = \sum_{a\le \sqrt{x}} \chi(a) \sum_{b\le x/a} \overline{\chi(b)} 
+\sum_{b\le \sqrt{x}} \overline{\chi(b)} \sum_{a\le x/b} \chi(a) - 
\sum_{a,b\le \sqrt{x}} \chi(a)\overline{\chi(b)}, 
$$
and (6) follows upon invoking the P{\' o}lya-Vinogradov bound (4).

Now we write $d(n)=\sum_{\ell | n} d_{\chi}(n/\ell) h(\ell)$ where 
$h$ is a multiplicative function with $h(p)=2-2\text{Re }\chi(p)$, and 
$|h(n)|\le d_4(n)$ for all $n$.  Observe that 
$$
x\log x+O(x)= \sum_{n\le x} d(n) 
= \sum_{\ell \le x} h(\ell) \sum_{m\le x/\ell} d_{\chi}(m).
$$
When $\ell \le x/q^2$ we use (6) to estimate the sum over $m$.  When $\ell$ is 
larger we trivially bound the sum over $m$ by $(x/\ell) \log (x/\ell) +O(x/\ell)$.   Thus 
we deduce that 
$$
x\log x+ O(x) \ll  \sum_{\ell \le x/q^2} |h(\ell)| \sqrt{xq/\ell}\log q + 
\sum_{x/q^2 \le \ell \le x} |h(\ell)| \frac{x}{\ell} \log q
\ll x\log q \sum_{\ell \le x} \frac{|h(\ell)|}{\ell}.
$$
 Since $\sum_{n\le x} |h(\ell)|/\ell \ll \exp(\sum_{p\le x} |h(p)|/p)
 =\exp(2{\Bbb D}(1,\chi;x)^2)$ we obtain the first part of the Lemma.

 To deduce the second part, note that the triangle inequality gives 
 $$
 ({\Bbb D}(f,\chi;x)+ {\Bbb D}(f,\psi;x))^2
 \ge \sum_{p\le x} \frac{1-\text{Re }|f(p)|^2 \chi(p)\overline{\psi}(p)}{p} 
 \ge \frac{1}{2}\sum_{p\le x} \frac{1-\text{Re }\eta(p)}{p},
 $$
 where $\eta$ is the primitive character of conductor below $Q^2$ which 
 induces $\chi \overline{\psi}$.  Now we appeal to the first part 
 of the   Lemma.
\enddemo  

\proclaim{Proposition 7} Let $\chi \pmod q$ be a primitive character and $t\in \Bbb R$.  
There is an absolute  constant $c>0$ such that for all 
 $x\ge q$ we have
 $$
 {\Bbb D}(1,\chi(n)n^{it};x)^2 \ge \frac12 \log \Big( \frac{c\log x}{\log (q(1+|t|))} \Big) .
 $$
Consequently, if $f$ is a multiplicative function, and $\chi$ and $\psi$ are any 
two distinct primitive characters with conductor below $Q$,  then for $x \ge Q$ we 
have
$$
{\Bbb D}(f,\chi(n)n^{it};x)^2 + {\Bbb D}(f,\psi(n) n^{iu};x)^2 \ge \frac{1}{8} \log \left( \frac{c\log x}{2\log (Q(1+|t-u|))} \right) .
$$
\endproclaim
\demo{Proof} The proof is much like that of Proposition 6, with some small changes. 
In place of $d_{\chi}(n)$ we will consider $d_{\chi,t}(n)=\sum_{ab=n} \chi(a)a^{it}
\overline{\chi(b)}b^{-it}$, and require an estimate like (6).  
To do this, we note that partial summation and the P{\' o}lya-Vinogradov inequality 
(4) yield 
$$
\align
\sum_{n\leq x} \chi(n)n^{it}&= x^{it} \sum_{n\leq x} \chi(n)
- it \int_1^x u^{it-1} \sum_{n\leq u} \chi(n) du 
\ll \sqrt{q}\log q  ( 1+ |t|\log x).\\
\endalign
$$
Using this, and arguing as in (6), we obtain 
$$
\sum_{n\le x} d_{\chi,t}(n) \ll \sqrt{qx} \log q (1+|t|\log x) + q \log^2 q (1+|t|\log x)^2. 
$$
The rest of the proof follows the lines of Proposition 6, breaking now into the 
cases when $\ell \le x/(q^2(1+|t|)^2)$, and when $\ell$ is larger.
\enddemo

\enddocument